
\input amstex.tex
\documentstyle{amsppt}
\magnification=\magstep1
\hsize=12.5cm
\vsize=18cm
\hoffset=1cm
\voffset=2cm

\footline={\hss{\vbox to 2cm{\vfil\hbox{\rm\folio}}}\hss}
\nopagenumbers
\def\DJ{\leavevmode\setbox0=\hbox{D}\kern0pt\rlap
{\kern.04em\raise.188\ht0\hbox{-}}D}

\def\txt#1{{\textstyle{#1}}}
\baselineskip=13pt
\def\hf{{\textstyle{1\over2}}}
\def\a{\alpha}
\def\d{{\,\roman d}}
\def\e{\varepsilon}

\def\={\;=\;}

\def\z{\zeta}

\def\hf{{\textstyle{1\over2}}}
\def\txt#1{{\textstyle{#1}}}

\font\tenmsb=msbm10
\font\sevenmsb=msbm7
\font\fivemsb=msbm5
\newfam\msbfam
\textfont\msbfam=\tenmsb
\scriptfont\msbfam=\sevenmsb
\scriptscriptfont\msbfam=\fivemsb
\def\Bbb#1{{\fam\msbfam #1}}

\def \NN {\Bbb N}

\def \RR {\Bbb R}

\font\ff=cmr8
\def\txt#1{{\textstyle{#1}}}
\baselineskip=13pt

\font\teneufm=eufm10
\font\seveneufm=eufm7
\font\fiveeufm=eufm5
\newfam\eufmfam
\textfont\eufmfam=\teneufm
\scriptfont\eufmfam=\seveneufm
\scriptscriptfont\eufmfam=\fiveeufm
\def\mathfrak#1{{\fam\eufmfam\relax#1}}

\font\tenmsb=msbm10
\font\sevenmsb=msbm7
\font\fivemsb=msbm5
\newfam\msbfam
     \textfont\msbfam=\tenmsb
      \scriptfont\msbfam=\sevenmsb
      \scriptscriptfont\msbfam=\fivemsb
\def\Bbb#1{{\fam\msbfam #1}}

\def \NN {\Bbb N}

\def \RR {\Bbb R}

  \def\rightheadline{{\hfil{\ff
  On a problem of Erd\H os}\hfil\tenrm\folio}}

  \def\leftheadline{{\tenrm\folio\hfil{\ff
   A. Ivi\'c }\hfil}}
  \def\emptyheadline{\hfil}
  \headline{\ifnum\pageno=1 \emptyheadline\else
  \ifodd\pageno \rightheadline \else \leftheadline\fi\fi}

\topmatter
\title
ON A PROBLEM OF ERD\H OS INVOLVING THE LARGEST PRIME FACTOR OF $n$
\endtitle
\author   Aleksandar Ivi\'c  \endauthor
\address
Aleksandar Ivi\'c, Katedra Matematike RGF-a
Universiteta u Beogradu, \DJ u\v sina 7, 11000 Beograd,
Serbia (Yugoslavia).
\endaddress
\email {\tt
aivic\@rgf.bg.ac.yu,  aivic\@matf.bg.ac.yu} \endemail
\dedicatory
\enddedicatory
\abstract {Let $P(n)$ denote the largest prime factor of an
integer $n \;(\ge 2)$, and let $N(x)$ denote the number of natural
numbers $n$ such that $2 \le n \le x$, and $n$ does not divide $P(n)!$.
We prove that $$ N(x) = x\left(2+O\left(\sqrt{\log_2x/\log x}\,\right)
\right)\int_2^x\rho(\log x/\log t) {\log
t\over t^2}\d t, $$ where $\rho(u)$ is the Dickman-de Bruijn
function. In terms of elementary functions we have $$ N(x) =
x\exp\left\{-\sqrt{2\log x\log_2 x}\left(1 +
O(\log_3x/\log_2x)\right)\right\}, $$ thereby sharpening  and
correcting recent results of  K. Ford and J.-M. De Koninck and N.
Doyon. }
\endabstract
\endtopmatter

\head
1. Introduction and statement of results
\endhead
Let $P(n)$ denote the largest prime factor of an integer $n \;(\ge 2)$.
In 1991 P. Erd\H os [5] proposed the following problem: prove that
$N(x) = o(x)\;(x\to\infty)$, where $N(x)$ denotes the number of
natural numbers $n$ such that $2 \le n \le x$, and $n$ does not divide
$P(n)!$. This problem is connected to the so-called Smarandache function
$S(n)$, the smallest integer $k$ such that $n|k!$.

Erd\H os's assertion was shown to be true by I. Kastanas [13],
and S. Akbik [1] proved later that
$N(x) \ll x\exp(-{1\over4}\sqrt{\log x}\,)$
holds. K. Ford [8] proposed an asymptotic formula for $N(x)$.
His Theorem 1 states that
$$
N(x) \sim {{\sqrt\pi}(1+\log2)\over 2^{3/4}}(\log x\log_2x)^{3/4}
x^{1-1/u_0}\rho(u_0)\qquad(x\to\infty).\leqno(1.1)
$$
In this formula one should correct the constant, as will be show in
Section 4. Here, as usual,
$\log_1 x \equiv \log x$ and $\log_k x = \log(\log_{k-1}x)$ for $k\ge 2$.
The function $\rho(u)$ is the continuous solution to the difference delay
equation $u\rho'(u) = -\rho(u-1)$ with the initial condition $\rho(u) = 1$
for $0 \le u \le 1$ and is commonly called the Dickman (or
Dickman-de Bruijn) function. The function $u_0 = u_0(x)$ is
implicitly defined by the equation
$$
\log x = u_0\left(x^{1/u_0^2}-1\right),
$$
so that
$$
u_0(x) = \left({2\log x\over\log_2 x}\right)^{1/2}
\left\{1 - {\log_3x\over2\log_2x} + {\log2\over2\log_2x}
+ O\left(\left({\log_3x\over\log_2x}\right)^2\right)\right\}.\leqno(1.2)
$$
The asymptotic formula (1.1) and the explicit expression (1.2) enabled
Ford (op. cit.) in Corollary 2 to deduce that
$$
N(x) = x\exp\left\{-(\sqrt{2}+o(1))\sqrt{\log x\log_2x}\right\}
\qquad(x\to\infty).\leqno(1.3)
$$

Unaware of Ford's work, J.-M. De Koninck and N. Doyon [3] recently
tackled this problem again. In [3] they published the result that
$$
N(x) = x\exp\left\{-(2+o(1))\sqrt{\log x\log_2x}\right\}\qquad(x\to\infty).
\leqno(1.4)
$$
Unfortunately, (1.4) is not true (it contradicts (1.3)).
Namely the argument in [3] leading to their (3.8) does not
hold, and the crucial parameter $u = \log x/\log y$ in the formula for
$\Psi(x,y)$ has to be evaluated more carefully. A straightforward
modification of their proof leads then again to (1.3).

\bigskip
The aim of this note is to  sharpen (1.1) and (1.3), and to provide
asymptotic formulas for sums of $S^r(n)$ and $1/S^r(n)$ when
$r\in\RR$ is fixed.
The results are contained in  the following

\bigskip
THEOREM 1.
$$
N(x) = x\exp\left\{-\sqrt{2}L(x)\left(1 + g_0(x)
+ O\left(\left({\log_3x\over\log_2x}\right)^3\right)\right)\right\},
\leqno(1.5)
$$
{\it where}
$$
L \;=\;L(x) \;=\;\sqrt{\log x\log\log x},\leqno(1.6)
$$
{\it and, for $r > -1$},
$$\eqalign{
g_r(x) &= {\log_3x+\log(1+r)-2-\log 2\over2\log_2x}
\left(1+{2\over\log_3x}\right) \cr&-
{(\log_3x + \log(1+r)-\log2)^2\over8\log_2^2x},\cr}\leqno(1.7)
$$

\bigskip
THEOREM 2. $$ N(x) = x\left(2 +
O\left(\sqrt{\log_2x\over\log x}\,\right)\right)
\int_2^x\rho\left({\log x\over\log t}\right) {\log t\over t^2}\d
t.\leqno(1.8) $$

\bigskip
THEOREM 3.
{\it Let $r>0$ be a fixed number. Then we have}
$$
\sum_{2\le n\le x}\,{1\over S^r(n)} =
x\exp\left\{-\sqrt{2r}L(x)\left(1+g_{r-1}(x) + O\left(\left(
{\log_3x\over\log_2x}\right)^2\right)\right)\right\},\leqno(1.9)
$$
{\it where $L(x)$ and $g_r(x)$ are given by} (1.6) {\it and} (1.7),
{\it respectively. Moreover, for any fixed integer $J \ge 1$
there exist computable constants $a_{1,r} = \z(r+1)/(r+1)$, \
$a_{2,r}\ldots,a_{J,r}$ such that}
$$
\sum_{2\le n\le x} S^r(n) =
x^{r+1}\sum_{j=1}^J\,{a_{j,r}\over\log^j x} +
O\left({x^{r+1}\over\log^{J+1}x}\right).\leqno(1.10)
$$

The asymptotic formula (1.8) is sharper than (1.5)--(1.7), since
it is a true asymptotic formula, while (1.5)--(1.7) is actually an
asymptotic formula for $\log {N(x)\over x}$. On the other hand,
the right-hand side of (1.5) is given in terms of elementary
functions, while (1.8) contains the (non-elementary) function
$\rho(u)$. Thus is seemed in place to give both (1.5) and (1.8),
especially since the proof of (1.5) requires only the more
elementary analysis which follows the work of De Koninck--Doyon
[3].

The asymptotic formula (1.9) seems to be new, while (1.10) sharpens
the results of S. Finch [7], who obtained (1.10) in the case $J=1$.

\medskip
{\bf Acknowledgement.} I thank Prof. G. Tenenbaum for
pointing out formula (2.6) to me, and indicating how it
improves the error term $O(\log_3x/\log_2x)$ in  the first version
of Th. 2 (coming from (2.5)), to $O(\sqrt{\log_2x/\log x})$.
I also thank the referee for valuable remarks.

\head
2. The necessary results of $\Psi(x,y)$
\endhead

The proofs of both Theorem 1 and Theorem 2 depend on results
on the function
$$
\Psi(x,y) \;=\;\sum_{n\le x ,P(n)\le y}\,1\qquad(2 \le y \le x),\leqno(2.1)
$$
which represents the number of $n$ not exceeding $x$, all of whose
prime factors do not exceed $y$. All the results that follow  can be
found e.g., in [9], [10] and  [14]. We have the elementary bound
$$
\Psi(x,y) \ll x\exp\left(-{\log x\over2\log y}\right)\qquad(2\le y \le x)
\leqno(2.2)
$$
and the more refined asymptotic formula
$$\eqalign{
\Psi(x,y) &= x\rho(u)\left(1 + O\left({\log(u+1)\over\log y}\right)\right)
 \cr u &= {\log x/\log y},\;\exp((\log_2x)^{5/3+\e}) \le y \le x.\cr}
\leqno(2.3) $$ Note that the Dickman--de Bruijn function $\rho(u)$
admits an asymptotic expansion, as $u\to\infty$, which in a
simplified form reads $$ \rho(u) = \exp\left\{-u\left(\log u +
\log_2u - 1 + O\left({\log_2u\over\log u}\right)\right)\right\}.
$$ Let, for given $u>1$, $\xi = \xi(u)$ be defined by $u\xi =
{\roman e}^\xi -1$, so that $$ \xi(u) = \log u +\log_2u +
{\log_2u\over\log u} + O\left(\left({\log_2u\over\log
u}\right)^2\right).\leqno(2.4) $$ If $u >2,\,|v| \le \hf u$, then
we have the asymptotic formula $$ \rho(u-v) =
\rho(u)\exp\left\{v\xi(u) + O((1+v^2)/u)\right\}.\leqno(2.5) $$
This result (see e.g., Hildebrand--Tenenbaum [10])
for bounded $|v|$ (see (3.3.9) of de la
Bret\`eche--Tenenbaum [2]) yields the asymptotic formula
$$ \rho(u-v)
=\rho(u){\roman e}^{v\xi(u)}\cdot\left\{1+O\left({1 \over
u}\right)\right\}\qquad(|v| \le v_0).\leqno(2.6)
$$

\head 3.
Proof of Theorem 1
\endhead
We turn now first to the proof of Theorem 1. The idea is to obtain
upper and lower bounds of the form given by (1.5)--(1.7).
 For the former,
note that if $P^2(n)|n$, then $n$ is counted by $N(x)$, since $P^2(n)$
cannot divide $P(n)!$. Therefore, if $T_0(x)$ denotes the number of
natural numbers $n\le x$ such that $P^2(n)|n$, then
$N(x) \ge T_0(x)$. The desired lower bound follows then from the
formula for $T_0(x)$, which is a special case of the formula for
($r\ge 0$ is a given constant)
$$
T_r(x) \;:=\;\sum_{2\le n\le x,P^2(n)|n}\,{1\over P^r(n)},\leqno(3.1)
$$
obtained by C. Pomerance and the author [12], and sharpened by the author
in [11]. Namely for $T_r(x)$ it was shown that
$$
T_r(x) = x\exp\left\{-(2r+2)^{1/2}L(x)\left(1+g_r(x) + O\left(\left(
{\log_3x\over\log_2x}\right)^3\right)\right)\right\},\leqno(3.2)
$$
where $L(x)$ and $g_r(x)$ are given by (1.6) and (1.7).
We note that Erd\H os,
Pomerance and the author in [6] investigated the related problem of
the sum of reciprocals of $P(n)$. They proved that (this result
is also sharpened in [11])
$$
\sum_{2\le n\le x}\,{1\over P(n)} = \sum_{p\le x}{1\over p}\Psi\left(
{x\over p},p\right) =
x\left(1 + O\left(\sqrt{\log_2x\over\log x}\,\right)\right)
\int_2^x\rho\left({\log x\over\log t}\right){1\over t^2}\d t.\leqno(3.3)
$$
The same argument leads without difficulty to
$$
T_0(x) = \sum_{p\le x}\Psi\left({x\over p^2},p\right) =
x\left(1 + O\left(\sqrt{\log_2x\over\log x}\,\right)\right)
\int_2^x\rho\left({\log x\over\log t}\right){\log t\over t^2}\d t.\leqno(3.4)
$$
\medskip
It remains yet to deal with the upper bound for $N(x)$. As in [3], we
shall use the inequality ($p$ denotes primes)
$$
N(x) \le \sum_{2\le r \le \log x/\log2}\,\sum_{p\le x^{1/r}}\Psi\left(
{x\over p^r},pr\right).\leqno(3.5)
$$
To see that (3.5) holds note that, if $n$ is counted by $N(x)$, then
$n$ must be divisible by $p^r$ which does not divide $P(n)!$. The
condition $r\ge 2$ is necessary, since squarefree numbers $q\;(>1)$
divide $P(q)!$. Moreover, $pr\le P(n)$ cannot hold, since otherwise
the numbers $p,2p,\ldots,rp$ would all divide $P(n)!$, and so would
$p^r$, which is contrary to our assumption.
Thus $n = p^rm\,(r\ge2)$, $P(m) \le P(n) < pr$, which easily
gives (3.5).

\medskip
We restrict now the ranges of $r$ and $p$ on the right-hand of (3.5)
to the ones which will yield the largest contribution.

The contribution of $r > 3L$ is trivially
$$
\ll \sum_{r>3L}\,\sum_p xp^{-r} \ll x2^{-3L} = x\exp(-3\log2\cdot L),
$$
which is negligible, since $3\log 2 > \sqrt{2}$. Similarly,
the contribution of $p > \exp(2L)$ is negligible, and also for
$p > \exp((\log x)^{1/3})$ and $r > (\log x)^{1/6}\log_2x$ the
contribution is, for some $C>0$,
$$
\ll \sum_{r > (\log x)^{1/6}\log_2x}\sum_{p > \exp((\log x)^{1/3})}
xp^{-r} \ll x\exp(-C\sqrt{\log x}\cdot\log_2x),
$$
which is negligible.

With (3.5) and (2.3) we see that in the range
$$
L_1 := \exp((\log x)^{1/3}) < p \le \exp({\txt{1\over4}}L),
\qquad r \le (\log x)^{1/6}\log_2x
$$
we have
$$
\eqalign{
u & = {\log(x/p^r)\over\log pr} = {\log x - r\log p\over \log p + \log r}
\cr&
= {\log x + O((\log x)^{2/3+\e})\over\log p + O(\log_2x)}
= {\log x\over\log p}\left(1 +
O\left((\log x)^{\e-1/3}\right)\right).\cr}
$$
Therefore by (2.3) we obtain, as $x\to\infty$,
$$
\Psi\left({x\over p^r},pr\right) \le
{x\over p^r}\exp\left(-(1+o(1)){\log x\over \log p}
\log\left({\log x\over\log p}\right)\right),
$$
which gives
$$\eqalign{&
\sum_{L_1<p\le\exp({1\over4}L)}\,
\sum_{r\ge2}\Psi\left({x\over p^r},pr\right)\cr&
\ll x\sum_{L_1<p\le\exp({1\over4}L)}
{1\over p^2}\max_{L_1<p\le\exp({1\over4}L)}
\exp\left(-(1+o(1)){\log x\over\log p}
\log\left({\log x\over\log p}\right)\right)\cr&
\ll x\exp(-(2+o(1))L),\cr}
$$
which is negligible. The contribution of $p \le L_1, r >
(\log x)^{1/6}\log_2x$ is easily seen to be also negligible.
This means that the main contribution to the
right-hand side of (3.5) comes from the range
$$
\exp({\txt{1\over4}}L) \le p \le \exp(2L),\qquad 2 \le r \le
(\log x)^{1/6}\log_2x,\leqno(3.6)
$$
or more precisely, for some $C >\sqrt{2}$,
$$
\eqalign{
N(x) &\ll x\exp(-CL) 
\cr& + \sum_{2 \le r \le (\log x)^{1/6}\log_2x}\,
\sum_{\exp({\txt{1\over4}}L) \le p \le \exp(2L)}
\Psi\left({x\over p^r},pr\right).\cr}\leqno(3.7)
$$
By trivial estimation it transpires that the contribution of $r\ge8$
in (3.7) will be negligible, so that we may write
$$
N(x) \ll x\exp(-CL) + L\max_{c>0}\sum_{2\le r\le7}{\Cal T}_{c,r}(x),\leqno(3.8)
$$
where, for a given constant $c>0$ and $r\in \NN$,
$$
{\Cal T}_{c,r}(x) \;:=\;\sum_{\exp(cL-1)<p\le\exp(cL)}
\Psi\left({x\over p^r},pr\right),
$$
since (again by (2.3)) the contribution of $c \ge c_0$, $c_0$ a large
positive constant, is negligible.
If $p$ is in the range indicated by ${\Cal T}_{c,r}(x)$, then
$$
\eqalign{
u &= {\log xp^{-r}\over\log pr} = {\log x\over cL}\left(1 +
O\left(\sqrt{{\log_2x\over\log x}}\,\right)\right),\cr
\log u &= \log_2x - \log c - \log L +
O\left(\sqrt{{\log_2x\over\log x}}\,\right).\cr}
$$
Therefore, for $2 \le r \le 7$,
$$
\eqalign{
{\Cal T}_{c,r}(x) &\ll \sum_{\exp(cL-1)<p\le\exp(cL)}{x\over p^r}
\exp\Bigl\{-{\log x\over cL}\Bigl(1 +
O\Bigl(\sqrt{{\log_2x\over\log x}}\,\Bigr)\Bigr)\times\cr&
\times (\hf\log_2x + \hf\log_3x + O(1))\Bigr\}\cr&
\ll x\exp\left\{-cL(r-1) - {1\over2c}L\left(1 +
O\left({\log_3x\over\log_2 x}\,\right)\right)\right\},\cr}
$$
so that the largest contribution is for $r=2$.
The function $c + 1/(2c)$ attains its minimal value $\sqrt{2}$
when $c = 1/\sqrt{2}$, hence
$$
{\Cal T}_{c,r}(x) \ll x\exp\left\{-\sqrt{2}L\left(1 +
O\left({\log_3x\over\log_2x}\right)\right)\right\},
$$
and we obtain
$$
N(x) \le x\exp\left\{-\sqrt{2\log x\log_2 x}\left(1 +
O\left({\log_3x\over\log_2x}\right)\right)\right\}.
\leqno(3.9)
$$
This is somewhat weaker than the result implied by Theorem 1. But we
saw that the main contribution to (3.7) comes from ${\Cal T}_{c,2}(x)$,
which is quite similar to the main sum appearing in the estimation
(lower bound) of $T_0(x)$ (cf. [12]). The only difference is that
in the case of ${\Cal T}_{c,2}(x)$ we have $\Psi(xp^{-2},2p)$, while
in the case of $T_0(x)$ the terms $\Psi(xp^{-2},p)$ appear.
In the relevant range for $p$ (i.e., $\exp({1\over4}L) \le p
\le  \exp(2L)$) this will not make much difference, so that
using the arguments [12]
we can obtain for $N(x)$  the upper bound implied by the
expression standing in (1.5)--(1.7). It is only for the sake of clearness
of exposition that we gave all the details for the slightly weaker upper
bound (3.9). This discussion completes the proof of Theorem 1.

\head
4. Proof of Theorem 2
\endhead
For the proof of Theorem 2 we need an expression that is sharper
than (3.5), which is only an upper bound. If $n$ is counted by $N(x)$
then

\smallskip
a) Either $P^2(n)$ divides $n$, and the contribution of such $n$ is
counted by $T_0(x)$ (cf. (3.1)), or

\smallskip
b) The number $n$ is of the form $n = pq^bk, p = P(n), b\ge 2,
P(k) < p, P(k) \not = q$, where henceforth $q$ (as well as $p$)
will denote primes. Since $n$ does not divide $P(n)!$, then
(similarly as in Section 3) it follows that $qb > P(n)$.  Also,
again as in Section 3, it is found that the contribution of $p \le
\exp(L/4) = Y_1$ is negligible. The contribution of $p > \exp(2L)
= Y_2$ is also negligible. Namely in this case, since $q^2|n, q >
P(n)/b \gg \exp(2L)/\log x$, the contribution is clearly
$$ \ll
x\sum_{p\gg\exp(2L)/\log x}p^{-2} \ll x\exp(-{\txt{3\over2}}L).
$$
Therefore, if b) holds, we need only consider
numbers $n$ for which there is a number $b\ge2$ and a prime $q
\in (p/b,p)$, such that $n = pq^bk, p = P(n), P(k) < p, P(k) \not =
q$. The contribution of $P(k) = q$ can be easily seen to be
negligible, as can be also shown for the one pertaining to $P(k) =
p$. Using the technique of Section 3 (or of [8]) it is seen
that the main contribution comes from the sum with $b=2$. If $b=2$, namely
if $n = pq^2k$, $p/2 < q < p$ and $q_1^2|k$ with some prime $p/2 < q_1 < p$,
then $n$ is counted at least twice, but it is seen that
the contribution of such $n$ is $\ll x\exp(-{\txt{3\over2}}L)$.
Therefore the contribution to $N(x)$ of integers satisfying b)
equals
$$
O(x\exp(-{\txt{3\over2}}L)) +
\sum_{Y_1<p\le Y_2} \sum_{p/2<q<p}\Psi\left({x\over
pq^2},p\right).\leqno(4.1)
$$
This is different from Ford [8, eq.
(2)], who had $\Psi\left({x\over pq^2},q\right)$ instead of
the correct $\Psi\left({x\over pq^2},p\right)$. This oversight does
affect his result (1.1), which is not correct since the constant
is not the right one ($1+\log 2$ should be replaced by 2),
 but the relation (1.3) remains valid, since both
expressions with the $\Psi$-function are of the same order of
magnitude. It follows that our starting relation for the proof
of Theorem 2 takes the shape
$$
N(x) = T_0(x) + O(x\exp(-{\txt{3\over2}}L)) +
\sum_{Y_1<p\le Y_2}\sum_{p/2<q<p} \Psi\left({x\over
pq^2},p\right).\leqno(4.2)
$$
By using the asymptotic formula
(2.3) for $\Psi(x,y)$, the prime number theorem in the standard
form (see e.g.,  [14])
$$
\pi(x) = \sum_{p\le x}1  = \int_2^x{\d t\over\log t} + O(x\exp(-\sqrt{\log
x})),\leqno(4.3)
$$
we write the sum in (4.2) as a Stieltjes
integral and integrate by parts. We set for brevity
$$ R :=
\sqrt{\log_2x\over\log x}.
$$
Then we see that the sum in question
equals
$$
\eqalign{& (1+O(R))x\int_{Y_1}^{Y_2}{1\over\log
t}\int_{t/2}^t {1\over ty^2\log y}\rho\left({\log x\over\log t} -
2{\log y\over\log t} -1 \right)\d y\d t\cr& =
(1+O(R))x\int_{Y_1}^{Y_2}{1\over t\log t}\int_{1-\log 2/\log t}^1
{1\over t^{z}z} \rho\left({\log x\over\log t} - 2z-1
\right)\d z\d t,\cr}
$$
on making the substitution $\log y = z\log
t$. Since both $z$ and $1/z$ equal $1 + O(R)$ in the relevant
range, by the use of (2.5) we see that our sum becomes
$$
\eqalign{& (1+O(R))x\int_{Y_1}^{Y_2}{{\roman e}^{2\xi(\log
x/\log t-1)}\over t\log t} \rho\left({\log x\over\log
t}-1\right)\int_{1-\log2/\log t}^1 t^{-z}\d z\d t \cr& =
x(1+O(R))\int_{Y_1}^{Y_2} {{\roman e}^{2\xi(\log
x/\log t-1)}\over t^2\log^2 t} \rho\left({\log x\over\log
t}-1\right)\d t.\cr}\leqno(4.4)
$$
Now we use (2.6) (with $u-3$
replacing $u$, $v=-2$), setting $u = \log x/\log t$, $Y_1 \le t
\le Y_2$, to obtain that $$\eqalign{& {\roman
e}^{2\xi(u-1)}\rho(u-1) = {\roman e}^{2\xi(u-1)}\rho((u-3)+2)\cr&
= {\roman e}^{2\xi(u-1)-2\xi(u-3)}\rho(u-3)\cdot(1 + O(R)).\cr}
$$
Since $\xi'(u) \sim {1/u}\;(u\to\infty)$ we have $$ {\roman
e}^{2\xi(u-1)-2\xi(u-3)}  = {\roman e}^{O(1/u)} = 1 + O\left(
{1\over u}\right) = 1 + O(R).
$$
Therefore the last integral in
(4.4) equals
$$ (1+O(R))\int_{Y_1}^{Y_2}\rho\left({\log x\over\log
t}-3\right){\d t\over t^2\log^2t}. \leqno(4.5)
$$
With the change of variable
$$
u = {\log x\over\log t},\; {\d t\over t\log^2t} = - {\d
u\over\log x}
$$
it follows that (4.5) becomes
$$
\eqalign{ S &:= {1\over\log x}(1+O(R))\int_{y_1}^{y_2}\rho(u-3)
\exp\left(-{\log x\over u}\right)\d u,\cr y_1 &:= \hf\sqrt{\log
x\over\log_2x},\quad y_2 := 4\sqrt{\log x\over\log_2x}.\cr}
$$
But since $u\rho'(u) = -\rho(u-1)$, integrating by parts we obtain
$$
\eqalign{ S &= {1\over\log
x}(1+O(R))\int_{y_1}^{y_2}-u\rho'(u-2) \exp\left(-{\log x\over
u}\right)\d u,\cr& = {1\over\log
x}(1+O(R))\int_{y_1}^{y_2}\rho(u-2) \left(1+{\log x\over
u}\right) \exp\left(-{\log x\over u}\right)\d u,\cr& =
(1+O(R))\int_{y_1}^{y_2}-(u-1)\rho'(u-1) \left({1
\over u}\right) \exp\left(-{\log x\over u}\right)\d u,\cr& =
(1+O(R))\int_{y_1}^{y_2}\rho(u-1) \left({\log x\over
u^2}\right) \exp\left(-{\log x\over u}\right)\d u,\cr&
= \log x(1+O(R))\int_{y_1}^{y_2}-\rho'(u)\left({1\over
u}\right) \exp\left(-{\log x\over u}\right)\d u,\cr& =
\log x(1+O(R))\int_{y_1}^{y_2}\rho(u) \left({\log x\over u^3}
\right) \exp\left(-{\log x\over u}\right)\d u,\cr}
$$
where we used several times that $u \asymp 1/R$ in the range of integration,
so that lower order terms could be absorbed by the $O(R)$--term.
Making again the change of variable $\log x/\log t = u$, we obtain
$$ S =
(1+O(R))\int_{Y_1}^{Y_2}\rho\left({\log x\over\log t}\right) {\log
t\over t^2} \d t = (1+ O(R)){T_0(x)\over x},
$$ where (3.4) was
used. This proves (1.8).

\vfill
\break
\bigskip
\head
5. Moments of the Smarandache function
\endhead
In this section we shall prove Theorem 3. We begin with the
proof of (1.9). Recall that the Smarandache
function $S(n)$ denotes the smallest $k\in \NN$ such that $n|k!$.
This implies that $P(n) \le S(n) \le n$, and if $S(n) \not = P(n)$,
then $n$ does not divide $P(n)!$. Thus we may write
$$
\eqalign{
\sum_{2\le n\le x}{1\over S^r(n)}& =
\sum_{2\le n\le x, S(n)=P(n)}{1\over S^r(n)}
+  \sum_{2\le n\le x,S(n)\not=P(n)}{1\over S^r(n)}\cr&
= \sum_{2\le n\le x}{1\over P^r(n)}
+ O\left(\sum_{2\le n\le x,n\not|P(n)!}{1\over P^r(n)}.\right).\cr}
\leqno(5.1)
$$
We have, by Ivi\'c--Pomerance [12],
$$
\sum_{2\le n\le x}{1\over P^r(n)} = x\exp\left\{
-\sqrt{2r}L(x)\Bigl(1+g_{r-1}(x)+O\Bigl(\Bigl({\log_3x\over\log_2x}\Bigr)^2
\Bigr)\Bigr)\right\},\leqno(5.2)
$$
with $L(x)$ and $g_r(x)$ given by (1.6) and (1.7), respectively.
On the other hand, by the argument that gives (3.5) and the proof
of Theorem 1,
$$
\eqalign{
&\sum_{2\le n\le x,n\not|P(n)!}{1\over P^r(n)}\cr&
\ll \sum_{2\le s\le\log x/\log2}\sum_{p\le x^{1/s}}{1\over p^r}
\Psi\left({x\over p^s},ps\right)\cr&
\ll x\exp\left\{-(\sqrt{2(r+1)}+o(1))L(x)\right\}
\qquad(x\to\infty).\cr}\leqno(5.3)
$$
The proof of (1.9) follows then from (5.1)--(5.3), with the remark
that a true asymptotic formula for the sum in (1.9) follows from
[11, eq. (2.11)], since the major contribution to the sum in
question comes from $n$ for which $S(n) = P(n)$. However, similarly
to the comparison between (1.5) and (1.8), the advantage of (1.9) is
that the right-hand side contains only elementary functions.

\smallskip
The proof of (1.10) utilizes the
following elementary lemmas, which generalize e.g.,
Lemma 1 and Lemma 2 of De Koninck--Ivi\'c [4].

\medskip
{\bf Lemma 1}. {\it Let $r>-1$ be a fixed number, and let
$c_{r,m} = (m-1)!(r+1)^{-m}$ for $m\in\NN$. For any fixed
integer $M\ge1$ we have}
$$
\sum_{p\le x}p^r = x^{r+1}\left({c_{r,1}\over\log x}
+ \ldots + {c_{r,M}\over\log^M x} + O\left({1\over\log^{M+1}x}\right)\right).
\leqno(5.4)
$$

\medskip
{\bf Proof}. Follows by partial summation from the prime number
theorem in the form (4.3).

\medskip
{\bf Lemma 2}. {\it For a given fixed real number $r>0$ and given
natural numbers $j$ and $m$, we have}
$$
\sum_{n\le{1\over2}x}{1\over n^{r+1}(\log x/n)^j}
= \sum_{k=0}^m\z^{(k)}(r+1){-j\choose k}{1\over(\log x)^{j+k}}
+ O\left({1\over\log^{j+m+1}x}\right).\leqno(5.5)
$$

\medskip
{\bf Proof}. Set $\ell = \ell(x) = \exp(\sqrt{\log x}).$ Then
$$
\sum_{n\le{1\over2}x}{1\over n^{r+1}(\log x/n)^j}
= \sum_{n\le\ell}{1\over n^{r+1}\log^jx
{\left(1 -{\log n\over\log x}\right)}^j} + O(\ell^{-r}).
$$
The assertion of the lemma follows if we use the binomial expansion
$$
{\left(1 -{\log n\over\log x}\right)}^{-j}
= \sum_{k=0}^m(-1)^k{ -j\choose k}
{\left({\log n\over\log x}\right)}^k + O\left(\left({\log n\over\log x}
\right)^{m+1}\right)
$$
and  ($\z(s)$ denotes the Riemann zeta-function)
$$
\sum_{n=1}^\infty \log^kn\cdot n^{-r-1} = (-1)^k\z^{(k)}(r+1)
\qquad(r > 0; k = 0,1,\ldots\,).
$$

\medskip
Now we turn to the proof of (1.10). We have, by Theorem 1,
$$\eqalign{
\sum_{2\le n\le x}S^r(n)&= \sum_{2\le n\le x,S(n)=P(n)}S^r(n)
+ \sum_{2\le n\le x,S(n)\not=P(n)}S^r(n)\cr&
= \sum_{2\le n\le x,S(n)=P(n)}P^r(n)
+ O\left(x^r\sum_{2\le n\le x,n\not|P(n)!}1\right)\cr&
= \sum_{2\le n\le x}P^r(n)
+ O\left\{x^{1+r}\exp(-(\sqrt{2}+o(1))L(x))\right\}\quad(x\to\infty).\cr}
\leqno(5.6)
$$
Observing that $\Psi(x/p,p) = [x/p]$ if $\sqrt{x}\le p\le x$, it follows
that ($k,n\in \NN, r>0$)
$$
\eqalign{
\sum_{2\le n\le x}P^r(n) &= \sum_{pk\le x,P(k)\le p}p^r
= \sum_{p\le x}p^r\Psi\left({x\over p},p\right)\cr&
=  \sum_{p\le \sqrt{x}}p^r\Psi\left({x\over p},p\right)
+ \sum_{\sqrt{x}<p\le x}p^r\Psi\left({x\over p},p\right)\cr&
= O\Bigl(x\sum_{p\le \sqrt{x}}p^{r-1}\Bigr)
+ \sum_{\sqrt{x}<p\le x}p^r\left[{x\over p}\right]\cr&
= O(x^{1+r/2}) + \sum_{pn\le x}p^r.\cr}\leqno(5.7)
$$
By using Lemma 1 and Lemma 2 we infer that
$$
\eqalign{&
\sum_{pn\le x}p^r =  \sum_{n\le{1\over2}x}\,\sum_{p\le x/n}p^r\cr&
= \sum_{n\le{1\over2}x} \left({x\over n}\right)^{r+1}
\left\{\sum_{j=1}^J {c_{r,j}\over\log^j(x/n)}
+ O\left({1\over\log^{J+1}(x/n)}\right)\right\} \cr&
= x^{r+1} \sum_{j=1}^J c_{r,j} \sum_{n\le{1\over2}x}
{1\over n^{r+1}\log^j(x/n)}  + O\left({x^{r+1}\over\log^{J+1}x}\right)\cr&
= x^{r+1} \sum_{j=1}^J c_{r,j}\sum_{i=0}^J\z^{(i)}(r+1) {-i\choose j}
{1\over(\log x)^{i+j}} + O\left({x^{r+1}\over\log^{J+1}x}\right) \cr&
= x^{r+1} \sum_{k=1}^J  {a_{k,r}\over\log^kx}
+ O\left({x^{r+1}\over\log^{J+1}x}\right)
 \cr}
$$
with
$$
a_{k,r} := \sum_{i=0}^{k-1}(k-i-1)!(r+1)^{i-k}\z^{(i)}(r+1){-i\choose k-i}.
\leqno(5.8)
$$
The assertion (1.10) of Theorem 3 follows then from (5.6)--(5.8), clearly one
has $a_{1,r} = \z(r+1)/(r+1)$.

\bigskip
Finally we remark that, for $r\in\RR$ fixed and $x\to\infty$,
$$
\sum_{2\le n\le x}{\left({S(n)\over P(n)}\right)}^r
= x + O\left\{x\exp(-(\sqrt{2}+o(1)))L(x)\right\}.\leqno(5.9)
$$
Namely, if $n = p_1^{\a_1}\ldots p_k^{\a_k}$ with $p_1 = P(n)$ is
the canonical decomposition of $n$, then $n$ divides
$$
\{(\a_1 + \a_2 + \ldots + \a_k)p_1\}! \le 
\left\{P(n)\left(\left[{\log n\over\log 2}\right]+1\right)\right\}!.
$$
Thus we have $P(n) \le S(n) \ll P(n)\log n$, and one obtains without
difficulty (5.9) by the method of proof of Theorem 3.

\vfill
\eject
\topglue2cm
\Refs
\bigskip

\item{[1]} S. Akbik, {\it On a density problem of Erd\H os}, Int. J. Math.
Sci. {\bf22}(1999), No. 3, 655-658.

\item{[2]} R. de la Bret\`eche and G. Tenenbaum, Entiers friables:
in\'egalit\'e de Tur\'an--Kubilius et applications, to appear, see
{\tt TKfriable.pdf} at {\tt http://www.iecn.u-nancy.fr/~tenenb/}

\item{[3]} J.-M. De Koninck and N. Doyon, {\it On a thin set of integers
involving the largest prime factor function},  Int. J. Math. Math.
Sci. (2003), No. {\bf19}, 1185-1192.

\item{[4]} J.-M. De Koninck and A. Ivi\'c, {\it The distribution of the
average prime factor of an integer}, Arch. Math. {\bf43}(1984), 37-43.

\item{[5]} P. Erd\H os, {\it Problem} 6674. Amer. Math. Monthly {\bf98}
(1991), 965.

\item{[6]} P. Erd\H os, A. Ivi\'c and C. Pomerance, {\it
On sums of reciprocals involving the largest prime factor
of an integer}, Glasnik Mat. {\bf21}(41) (1986), 283-300.

\item{[7]} S.R. Finch, {\it Moments of the Smarandache function}, Smarandache
Notions Journal {\bf11}(2000), No. 1-3, 140-141.

\item{[8]} K. Ford, {\it The normal behavior of the Smarandache
function}, Smarandache Notions Journal {\bf10}(1999), 81-86.

\item{[9]} A. Hildebrand, {\it On the number of positive integers $\le x$
and free of prime factors $>y$}, J. Number Theory {\bf22}(1986), No. 3,
289-307.

\item{[10]} A. Hildebrand and G. Tenenbaum, {\it Integers without
large prime factors}, J. Th\'eorie Nombres Bordeaux {\bf5}(1993),
411-484.

\item{[11]} A. Ivi\'c, {\it On sums involving reciprocals of the largest
prime factor of an integer II}, Acta Arith. {\bf71}(1995), 229-251.

\item{[12]} A. Ivi\'c and C. Pomerance, {\it Estimates for certain sums
involving the largest prime factor of an integer}, in ``Colloquia
Math. Soc. J. Bolyai {\bf34}. Topics in Analytic Number Theory", Budapest
1981, North-Holland, 769-789.

\item{[13]} I. Kastanas, {\it Solution to problem} 6674,
Amer. Math. Monthly {\bf98} (1994), 179.

\item{ [14]} G. Tenenbaum,   {\it Introduction to Analytic and
Probabilstic Number Theory}, Cambridge Studies in Advanced Mathematics
{\bf46}, Cambridge University Press, Cambridge, 1995.

\bigskip
\bigskip

Aleksandar Ivi\'c

Katedra Matematike RGF-a

Universitet u Beogradu

\DJ u\v sina 7, 11000 Beograd, Serbia

{\tt aivic\@rgf.bg.ac.yu}

\endRefs
\bigskip

\vfill


\bye